\documentclass{amsart}

\usepackage{hyperref}
\usepackage{amsrefs}
\usepackage{comment}
\usepackage{graphicx} 
\usepackage{float} 

\usepackage[dvipsnames]{xcolor}
\usepackage{tikz}
\usepackage{multirow}
\allowdisplaybreaks

\title{A New Proof of Fine's Identity using \\ Wildberger's Polynomial Formula}

\theoremstyle{definition}
\newtheorem{theorem}{Theorem}

\newcommand{\tri}{\mathbin{\overline{\nabla}}}

\newcommand{\m}{{\mathbf m}}
\newcommand{\C}{{\mathbf c}}
\newcommand{\T}{{\mathbf t}}

\newcommand{\Sb}{{\mathbf S}}

\newcommand{\Tb}{{\mathbf T}}

\newcommand{\Sm}{\mathcal{S}}
\newcommand{\Tm}{\mathcal{T}}  
\newcommand{\OB}{\big[}
\newcommand{\CB}{\big]}
\newcommand{\mC}{C_{\m}}
\newcommand{\mV}{V_{\m}}
\newcommand{\mE}{E_{\m}}
\newcommand{\mF}{F_{\m}}
\newcommand{\BigComment}[1]{}

\renewcommand{\k}{{\bf k}}

\providecommand{\U}[1]{\protect\rule{.1in}{.1in}}

\begin{document}

\author{Dean Rubine}
\email{DeanRubineMath@gmail.com}

\begin{abstract}
In 1959, N. J. Fine showed that the sum of the multinomial coefficients corresponding to the partitions of a natural number $n$ into $r$ parts is a binomial coefficient:
$$
 \sum_{\substack{k_1 + k_2 + k_3 + {}\ldots = r  \\   k_1 + 2k_2 + 3k_3 + {}\ldots = n }} \binom{r}{k_1, k_2, k_3, \ldots } 
 = \binom{n - 1}{r - 1}
$$
Fine gives a rather pithy proof, though we're still stuck on the part that says, ``We begin with an important though obvious remark.''

In 2025, Wildberger and Rubine gave the series solution to the general polynomial, derived from
a non-associative algebra of roofed, subdivided polygons they call \textit{subdigons}.
We generalize subdigons to \textit{tubdigons}, which include 2-gons, and count tubdigons of a given type
two ways: through a simple counting argument (backed up by the combinatorics literature) and by using Wildberger's polynomial formula to solve the polynomial implied by the multiset specification of tubdigons.
Comparing corresponding terms yields Fine's Identity.
\end{abstract}
\maketitle

\section{Subdigons and Polynomial Equations}
We recap Wildberger's derivation.
Wildberger and Rubine~\cite{Wildberger2025} define a \textbf{subdigon} $s$ of \textbf{type} $\m=[m_2,m_3,m_4,\ldots]$ to be a convex planar polygon with a distinguished side called its \textbf{roof}, that is subdivided by non-intersecting diagonals into $m_2$ triangles, $m_3$ quadrilaterals, $m_4$ pentagons and so on.
Necessarily only a finite number of $m_i$ are non-zero.
The null subdigon, denoted $|$, has two vertices, one edge, no faces and a type $\m=[\,]$. 
A subdigon of type $\m$ has counts of vertices, edges and faces,
\begin{align}
\mV &=  2+ m_2 +  2 m_3 
+ \ldots , 
\ \ 
\\
\mE &= 1 + 2 m_2 +  3m_3
+ \ldots , 
\ \ 
\\
\mF&= m_2 + m_3 
+ \ldots   
\label{eqn:VmEmFm}
\end{align}
satisfying Euler's polytope formula $ \mV - \mE + \mF = 1$.

We define $\Sm_{\m}$ as the finite multiset of subdigons of type $\m$; then the ongoing multiset $\Sm$ of all subdigons is layered by type: 
\begin{equation}
   \Sm \equiv \sum_{\m \ge 0} \Sm_{\m}  . \qquad \qquad
\end{equation}
We could use sets as no elements are repeated, but multisets conveniently have addition.

Subdigons are built up recursively through a family of $k$-ary operators $\tri_k$, $k \ge 2$, where 
$\tri_k(s_1, s_2, \ldots, s_k)$ consists of a central roofed $(k+1)$-gon, with attached subdigons $s_1, s_2, \ldots,\allowbreak s_k$ adjoined along their roofs in a counterclockwise fashion.
A subdigon $s$ is either $|$, or has a central roofed $(k+1)$-gon for some natural $k \ge 2$,
then necessarily of the form 
\begin{equation}
s=\tri_k(s_1,s_2,\ldots,s_k)
\end{equation}
for unique subdigons $s_1,s_2,\ldots,s_k$.
The equation for $\Sm$, the multiset of subdigons, follows:
\begin{equation} \label{eqn:Sk} 
\Sm =  \ \OB \  | \ \CB \ +  \tri_2(\Sm,\Sm) +  \tri_3(\Sm,\Sm,\Sm)  +  \tri_4(\Sm,\Sm,\Sm,\Sm)  + \ldots 
\end{equation}
where $\tri_k$ is extended to multisets of subdigons as:
\begin{equation}
    \tri_k(M_1, M_2, \ldots, M_k) = \OB \tri_k(s_1, s_2, \ldots, s_k) \! : \ s_1 \in M_1, s_2 \in M_2, \ldots, s_k \in M_k  \CB .
\end{equation}
The \textbf{accounting monomial} of a subdigon $s$ of type $\m$ is defined to be:
\begin{equation}
\psi(s) \equiv t_2^{m_2} t_3^{m_3} t_4^{m_4} \cdots \allowbreak \equiv \T^{\m} .
\end{equation}
We have $\psi(|)=1$, and since $\tri_k$ adds a $(k+1)$-gon:
\begin{equation}
 \psi(\tri_k(s_1, s_2, \ldots, s_k))=t_k \, \psi(s_1) \psi(s_2) \cdots \psi(s_k)
 \end{equation}
where we see the multiplications are essentially vector additions of the type vectors.
$\Psi(M)$ maps a multiset of subdigons to a polynomial, $\Psi(M) \equiv \sum_{s\in M} \psi(s)$,
so via linearity: 
\begin{equation}
\Psi(\tri_k(M_1, M_2, \ldots, M_k) ) = t_k \,  \Psi(M_2) \Psi(M_3) \cdots \Psi(M_k).
\end{equation}
The \textbf{hyper-Catalan number} $\mC$ counts the number of subdigons of type $\m$, $\mC=|\Sm_{\m}|$, and its multivariate generating series is called the \textbf{subdigon polyseries} $\Sb$: 
\begin{equation}
\Sb = \Sb[t_2,t_3, \ldots] \equiv \Psi(\Sm) = \sum_{\m \ge 0} \mC \T^{\m} .
\end{equation}
Applying $\Psi$ to equation (\ref{eqn:Sk}) tells us $\Sb$ is the solution to the general geometric polynomial equation:
\begin{theorem} (Wildberger's Soft Polynomial Formula)  \label{thm:spf}
The generating series for the hyper-Catalan numbers $\Sb[t_2, t_3, \ldots] = \sum_{\m \ge 0} \mC \T^{\m}$ satisfies:
$$ \Sb = 1 + t_2 \Sb^2 + t_3 \Sb^3 + t_4 \Sb^4 + \ldots $$
Further,
$\displaystyle 
x =\frac {c_0}{c_1} \Sb \ [\frac{c_0c_2}{c_1^2},\ \frac{c^2_0c_3}{c_1^3},\ \frac{c^3_0c_4}{c_1^4},\ \ldots] = \sum_{\m \ge 0} \mC \dfrac{ c_0^{1 + m_2 + 2m_3 + \ldots}} { c_1^{1 + 2m_2 + 3m_3+ \ldots}} c_2^{m_2} c_3^{m_3}  \cdots
$ satisfies
$$
c_0 - c_1 x + c_2 x^2 + c_3 x^3 + c_4 x^4 + \ldots = 0 .
$$
\end{theorem}
These are formal series solutions; no claim of convergence is being made.

Defining $\m! \equiv m_2! m_3! \cdots$, we can write
Erd\'elyi and Etherington's~\cite{Erdelyi1940} formula for $\mC$ as:
\begin{equation} \label{eqn:mC}
\mC = \dfrac{(\mE - 1)!}{(\mV-1)! \, \m!}
\end{equation}
It is not necessary to have this closed form for our proof of Fine's Identity.  We can write Wildberger's Soft Polynomial Formula as
\begin{equation}
     x = \sum_{\m \ge 0} \ \mC \  \dfrac{ c_0^{\mV -1} } { c_1^{\mE} } \, \C^{\m} .
\end{equation}

\section{Tubdigons}
We generalize subdigons to \textbf{tubdigons}, which are identical except that 2-gons are allowed.  We can't really draw 2-gons with straight edges, so we'll use at least one curved edge for every 2-gon. 

Subdigons have an obvious bijection with plane trees with no unary nodes.
By allowing 2-gons, tubdigons have a bijection that includes unary nodes.

\newcommand{\tR}{R_{m_1; \m}}
\newcommand{\tE}{E_{m_1; \m}}
\newcommand{\tF}{F_{m_1; \m}}
\newcommand{\tV}{V_{m_1; \m}}

We'll write tubdigon types as a vector or subscript with a semicolon: $[m_1; \m ]=[m_1; m_2, m_3, \ldots]$, where $m_1$ counts the number of 2-gons, and $\m$ is the type vector for a subdigon (index starting at 2) so $m_2$ counts triangles, $m_3$ quadrilaterals, etc.  A tubdigon of type $[m_1; \m ]$ has 
\begin{align}
\tV &=\mV \textrm{ vertices}, \\ \tE&=m_1+\mE \textrm{ edges, and } \\ \tF&=m_1 + \mF \textrm{ faces},    \label{eqn:VkEkFk}
\end{align}
the natural generalizations that start at $m_1$ instead of $m_2$.

Call $\tR \equiv R[m_1; \m]$ the number of tubdigons of type $[m_1; \m ]$, $R$ for Raney~\cite{Raney1960}.
It's straightforward to count those in terms of the subdigons $s$ of type $\m$.
We have $m_1$ extra edges to use to add 2-gons to $s$.  
Each extra edge can double up (or triple up, etc.) an existing edge of $s$.  
So this is the classic stars and bars problem; we're putting $m_1$ indistinguishable balls into $\mE$ distinguishable bins. There are $\mE-1$ walls between the bins; those walls can go in any of $m_1 + \mE -1$ places.
We can do this for each subdigon of type $\m$ so:
\begin{equation} \label{eqn:mR}
 \tR = \binom{m_1 + \mE - 1}{\mE - 1} \mC =  \binom{m_1 + \mE - 1}{m_1} \mC
\end{equation}
which with equation (\ref{eqn:mC}) expands out to the natural generalization of $\mC$,
\begin{equation}
\tR = \dfrac{(\tE - 1)!}{(\tV-1)! \, m_1 ! \, \m!}
\end{equation}

This formula appears in the combinatorics literature, perhaps first in Raney, who counts well-formed expressions of k-ary function applications (whose parse trees are general planes trees),
then in Tutte~\cite{Tutte1964}, who counts plane trees of a given type. 
Kreweras~\cites{Kreweras1972, KrewerasEarnshaw2005} obtains both formulas; $\tR$ from counting types of 
non-crossing partitions of a cyclic graph, and $\mC$ when he considers non-crossing partitions without singletons.  Of course $R[0; \m]=\mC$.

\section{The Tubdigon Equation}

We won't echo the entire subdigon definition in full for tubdigons, leaving that to the reader.
We've added the unary $\tri_1(r)$ to our family of panelling operators, so
the equation for $\Tm$, the multiset of tubdigons, is:
\begin{equation} \label{eqn:Tm} 
\Tm =  \ \OB \  | \ \CB \ +  \tri_1(\Tm) + \tri_2(\Tm,\Tm) +  \tri_3(\Tm,\Tm,\Tm)  +  \tri_4(\Tm,\Tm,\Tm,\Tm)  + \ldots 
\end{equation}
For a tubdigon $r$ of type $[m_1; \m ]$ we redefine $\psi(r) \equiv t_1^{m_1} \T^{\m}$ and $\Psi(M)$ for multisets of tubdigons, accordingly.  Applying the redefined $\Psi$,
\begin{equation}
\Tb = 1 + t_1 \Tb +  t_2 \Tb^2 + t_3 \Tb^3 + t_4 \Tb^4 + \ldots
\end{equation}
where $\Tb$ is the generating sum for $\tR$,
\begin{equation}
\Tb \equiv \Psi(\Tm)=\sum_{[m_1; \m ] \ge 0} \tR  \, t_1^{m_1} \, \T^{\m}
\end{equation}
\section{Applying the polynomial formula}
Alternatively, to solve
\begin{equation}
0 = 1 - (1- t_1) \Tb +  t_2 \Tb^2 + t_3 \Tb^3 + t_4 \Tb^4 + \ldots
\end{equation}
we apply Theorem \ref{thm:spf}, the soft polynomial formula from Wildberger

\begin{equation}
 \Tb = \sum_{\m \ge 0} \mC \frac{1}{(1- t_1)^{\mE}} \T^{\m} 
\end{equation}
$1/(1-t_1)$ is the geometric series.
\begin{equation}
 \Tb = \sum_{\m \ge 0} \mC \left( \sum_{i \ge 0} t_1^ {i} \right) ^{\mE}  \T^{\m} 
\end{equation}
Applying the multinomial theorem,
\begin{align}
\Tb &= \sum_{\m \ge 0} \mC \left( \sum_{\substack{j_i \ge 0 \\ \sum_i j_i = \mE }} 
\binom{\mE}{j_0, j_1, \ldots } \prod_{i  \ge 0} (t_1^i) ^{j_i} \right) \T^{\m} 
\\ \Tb &= \sum_{\m \ge 0} \mC  
\sum_{\substack{ \sum_i j_i = \mE }}  \!\!\!\!\! \binom{\mE}{j_0, j_1, \ldots } t_1^{ \sum_{i \ge 0} i j_i} \T^{\m} 
\end{align}
We  drop $j_i \ge 0$ from the notation, though it of course remains true.

We extract coefficients from the generating series, $\tR =  [t_1 ^{m_1} ; \T^{\m} ] \Tb $.
\begin{equation}
\tR =\mC \sum_{\substack{ \sum_{i  \ge 0}  j_i = \mE  \\   \sum_{i  \ge 0} i j_i = m_1}} \binom{\mE}{j_0, j_1, j_2, \ldots } 
\end{equation}
We need to reconcile starting at $j_0$  (because the geometric series starts at $i=0$) not $j_1$ like in Fine's Identity.
Let $j_i=k_{i+1}$ (we assume and omit $k_i \ge 0$ from the notation).
\begin{equation}
\tR =\mC \sum_{\substack{ \sum_{i  \ge 1}  k_{i} = \mE  \\   \sum_{i  \ge 1} i k_{i} = m_1 + \mE }} \binom{\mE}{k_1, k_2, \ldots } 
\end{equation}
Comparing that to equation (\ref{eqn:mR}) for $\tR$, and noting $\mC \ne 0$, we have:
\begin{equation}
 \sum_{\substack{ \sum_{i  \ge 1}  k_{i} = \mE  \\   \sum_{i  \ge 1} i k_{i} = m_1 + \mE }} \binom{\mE}{k_1, k_2, \ldots } 
 =
\binom{m_1 + \mE - 1}{ m_1}
.
\end{equation}

Setting $n=m_1+\mE, \ r=\mE$, we've derived Fine's Identity~\cites{Fine1959,Fine1988} using Wildberger's Polynomial Formula and the multinomial theorem.
\begin{theorem}[Fine's Identity]
$$
 \sum_{\substack{ \sum_{i  \ge 1}  k_{i} = r  \\   \sum_{i  \ge 1} i k_{i} = n }} \binom{r}{k_1, k_2, \ldots } 
 =
\binom{n - 1}{r - 1}
$$
\end{theorem}

\newcommand\as[2]{a_{#1}[#2]}
\newcommand\aj[1]{\as{j}{#1}}

\section{``An Important Though Obvious Remark''}

Let's review Fine's proof.
It begins with a lemma presented without proof, ``an important though obvious remark,'' as he states (renotated, and a typo fixed):

\begin{theorem}[Fine's Lemma]\label{thm:obvious}
    Let $S_j(q) = \sum_{k \ge 0} \aj{k} \, q^k$ for $j=1,2,\ldots$.  Then
    $$
        \prod_{j \ge 1} S_j(q^j) = \sum_{n \ge 0} q^n \ \sum_{\substack{   \sum_{i  \ge 1} i k_{i} = n }} \ \ \prod_{j \ge 1} \aj{k_j}
    $$
\end{theorem}
We have a family of series $S_j(q)$, each the generating series of a sequence $\aj{k}$.  We pass $q^j$ as the argument to the $j$th series, and multiply them all together.  The result is the generating series for the sequence whose $n$th term is the sum over partitions of $n$ of $\prod_j \aj{k_j}$, the product of the sequence elements corresponding to the partition.

Fine's clever lemma is an excellent starting point for the subject of sums of partitions.
Like Fine, we won't attempt a proof.
Let's write out a few terms of the product of the first four to see what's going on.  We'll assume $\aj{0}=1$; i.e. all the sequences start with 1, which is the common case, and less cluttered:
\begin{align}
{\textstyle 
\prod_{j=1}^4 } & S_j(q^j) = (1 + \as 1 1 q + \as 1 2 q^2  + \as 1 3  q^3 +\as 1 4 q^4 + \ldots) 
\\& \qquad\quad \cdot
(1 + \as 2 1 q^2 + \as 2 2 q^4  + \as 2 3 q^6 + \as 2 4 q^8 +\ldots) 
\nonumber \\& \qquad\quad \cdot 
(1 + \as 3 1  q^3 + \as 3 2  q^6  + \as 3 3  q^9 + \as 3 4 q^{12}+ \ldots) 
\nonumber \\& \qquad\quad\cdot 
(1 + \as 4 1  q^4 + \as 4 2 q^8  + \as 4 3  q^{12} + \as 4 4 q^{16} +\ldots) 
\nonumber \\
& =  1 
+ \as 1  1   \, q 
+ (  \as 1 2 + \as 2 1  ) \,  q^2
+ (\as 1 3 + \as 1 1 \as 2 1 + \as 3 1  )  \, q^3 
\nonumber \\&
 \quad + (\as 1 4 + \as 1 2 \as 2 1 + \as 2 2 + \as 1 1 \as 3 1 + \as 4 1 ) \, q^4
+ \ldots 
\nonumber 
\end{align}

Each term multiplying $q^n$ corresponds to the partitions of $n$; the $q^0$ term is the product of the constants, unity here.

For the identity in question (equation (22.17) in \cite{Fine1988}), Fine applies the lemma to $S_j(q)=\exp(tq)$, all identical series, with $\aj{k} = t^k/k!$ by the usual expansion.  

Focusing on the left side of Theorem \ref{thm:obvious}, the exponential argument becomes a geometric series, which we sum and then expand the exponential, apply the binomial theorem with a $-r$ power, and finally apply the identity $\binom{-r}{m}=(-1)^m \binom{r+m-1}{r-1}$ and the substitution $n=m+r$.
\begin{align}
\prod_{j \ge 1}& S_j(q^j) 
=\prod_{j \ge 1} \exp(tq^j)
= \exp\left(\sum_{j \ge 1} tq^j\right) 
= \exp\left( \dfrac{tq}{1-q} \right)  
= \sum_{r \ge 0} \dfrac{1}{r!} \left( \dfrac{tq}{1-q} \right)  ^r 
\\&
= \sum_{r \ge 0} \dfrac{t^r}{r!} \left( \dfrac{ 1-q }{q} \right) ^{-r} 
= \sum_{r \ge 0} \dfrac{t^r}{r!} (  q^{-1} -1 ) ^{-r} 
= \sum_{r \ge 0} \dfrac{t^r }{r!}  \sum_{ m \ge 0} \binom{-r}{m} (q^{-1})^{-r-m} (-1)^{m} 
\nonumber \\&
= \sum_{r \ge 0} \dfrac{t^r }{r!}  \sum_{ m \ge 0} \binom{r+m-1}{r-1} q^{r+m} 
= \sum_{r \ge 0} \dfrac{t^r }{r!} \sum_{ n \ge 0} \binom{n-1}{r-1} q^{n}  
\nonumber 
\end{align}
The right side of Theorem \ref{thm:obvious} is straightforward here. We transform the sum over all partitions into a double sum, over the number of parts in the partition $r$, and over partitions of $r$ parts:
\begin{align}
\sum_{n \ge 0}  & q^n \sum_{\substack{ \sum_{i  \ge 1} i k_{i} = n }} \prod_{\ j \ge 1} \aj{k_j}
= \sum_{n \ge 0}  q^n \sum_{\substack{ \ \sum_{i  \ge 1} i k_{i} = n }} \prod_{\ j \ge 1} \dfrac{ t^{k_j}}{k_j!}
= \sum_{n \ge 0}  q^n \sum_{\substack{ \sum_{i  \ge 1} i k_{i} = n }}   \dfrac{t^{\sum_{j \ge 0} k_j} }{\prod_{j \ge 1} k_j!}
\\&
= \sum_{n \ge 0}  q^n \sum_{r \ge 0} \sum_{\substack{ \sum_{i  \ge 1}  k_{i} = r \\   \sum_{i  \ge 1} i k_{i} = n }}   \dfrac{t^r }{\prod_{j \ge 1} k_j!}
= \sum_{r \ge 0} \dfrac{t^r}{r!}   \sum_{n \ge 0} \sum_{\substack{ \sum_{i  \ge 1}  k_{i} = r \\   \sum_{i  \ge 1} i k_{i} = n }}   \dfrac{r! }{\prod_{j \ge 1} k_j!} \ q^n 
\nonumber
\end{align} 
Equating coefficients and recognizing the multinomial coefficient yields Fine's Identity:
\begin{equation}
 \binom{n-1}{r-1} =  \ \sum_{\substack{ \sum_{i  \ge 1}  k_{i} = r \\   \sum_{i  \ge 1} i k_{i} = n }}  \binom{r}{k_1, k_2, \ldots} \quad\checkmark
\end{equation}

Fine further comments that if we sum over $r$, summing a row of Pascal's Triangle, we find the sum of multinomial coefficients corresponding to all the partitions of $n$ is $2^{n-1}$:
\begin{equation}
\sum_{\substack{ \sum_{i  \ge 1} i k_{i} = n }}  \binom{k_1 + k_2 + k_3 + \ldots }{k_1, k_2, k_3, \ldots} = 2^{n-1}
\end{equation}

\newcommand{\Mb}{{\bf M}}
\renewcommand{\U}{{\bf u}}
\newcommand{\uk}{\U^{\k}}

\section{Partitions as Edge Layers of Multinomial Series}
Fine's lemma generates partitions via a product of an unending list of univariate series. 
As an alternative,
let's note a relation between partitions of a natural number and tubdigon types.
We relax our notation to remove the semicolon; our vectors here start at index one, we have tubdigon type $\k \equiv [k_1, k_2, \ldots]$ and variables $\U \equiv [u_1, u_2, \ldots]$, so that $\U^{\k} \equiv u_1^{k_1} u_2^{k_2} \cdots$, a monomial.  Each such monomial $\U^{\k}$ may be seen as a partition of $n=1k_1 + 2k_2 + 3k_3 + \ldots$ with $r=k_1 + k_2 + k_3 + \ldots$ parts.

The partitions of $n$ correspond to the tubdigon types with $n+1$ edges (because $E_{\k}=n+1$, equations \eqref{eqn:VmEmFm} and \eqref{eqn:VkEkFk}).
We group the terms of a multinomial with all possible monomials into edge layers.
We use
\begin{equation}
\Mb[u_1, u_2, \ldots]=\sum_{ \k \ge 0} \U^{\k}
   \end{equation}
a series with an unbounded number of variables, where every possible monomial occurs once with a coefficient of unity.
We can layer the series by the number of edges $E_{\k} = 1 + k_1  + 2 k_2 + 3 k_3 + \ldots $ \cite{Rubine2025Finite}, by multiplying each $u_i$ by $e^i$, ($e$ being a literal variable whose exponent counts edges less one, and $i$ being the index), then collecting terms in powers of $e$, giving:
\begin{align}
\Mb[e^1 u_1 ,& \, e^2 u_2,  \,e^3 u_3, \, \ldots] = 1 
\nonumber  
+e^1
( u_{1})
+e^{2}
( u_{1}^{2} + u_{2})
+e^{3}
( u_{1}^{3} + u_{1} u_{2} + u_{3} )
\nonumber \\[0pt] {}
+e^{4} &( u_{1}^{4} + u_{1}^{2} u_{2} + u_{1} u_{3} + u_{2}^{2} + u_{4})
\nonumber \\[0pt] {}
+e^{5} &( u_{1}^{5} + u_{1}^{3} u_{2} + u_{1}^{2} u_{3} + u_{1} u_{2}^{2} + u_{1} u_{4} + u_{2} u_{3} + u_{5})
\nonumber \\[0pt] {}
+e^{6} & (u_{1}^{6} + u_{1}^{4} u_{2} + u_{1}^{3} u_{3} + u_{1}^{2} u_{2}^{2} + u_{1}^{2} u_{4} + u_{1} u_{2} u_{3} + u_{1} u_{5} + u_{2}^{3} + u_{2} u_{4} + u_{3}^{2} + u_{6} )
\nonumber \\[0pt] {}
+e^{7} &( u_{1}^{7} + u_{1}^{5} u_{2} + u_{1}^{4} u_{3} + u_{1}^{3} u_{2}^{2} + u_{1}^{3} u_{4} + u_{1}^{2} u_{2} u_{3} + u_{1}^{2} u_{5} 
\nonumber \\[0pt]  {} & \quad
+ u_{1} u_{2}^{3} + u_{1} u_{2} u_{4} + u_{1} u_{3}^{2} + u_{1} u_{6} + u_{2}^{2} u_{3} + u_{2} u_{5} + u_{3} u_{4} + u_{7})
\nonumber \\[0pt] {}
+e^{8} &( u_{1}^{8} + u_{1}^{6} u_{2} + u_{1}^{5} u_{3} + u_{1}^{4} u_{2}^{2} + u_{1}^{4} u_{4} + u_{1}^{3} u_{2} u_{3} + u_{1}^{3} u_{5} 
\nonumber \\[0pt] {} &\quad
+ u_{1}^{2} u_{2}^{3} + u_{1}^{2} u_{2} u_{4} + u_{1}^{2} u_{3}^{2} + u_{1}^{2} u_{6} + u_{1} u_{2}^{2} u_{3} + u_{1} u_{2} u_{5} + u_{1} u_{3} u_{4} 
\nonumber \\[0pt] &\quad
+ u_{1} u_{7} + u_{2}^{4} + u_{2}^{2} u_{4} + u_{2} u_{3}^{2} + u_{2} u_{6} + u_{3} u_{5} + u_{4}^{2} + u_{8}) + \cdots
\end{align}
where the terms multiplying $e^n$ are the partitions of $n$, e.g the $e^5$ row says $1+1+1+1+1=1+1+1+2=1+1+3=1+4=2+3=5$.

The number of terms at each level, gotten by setting the original variables $u_i$ to one, generates the partition function (A41 in OEIS):
\begin{align}
\Mb [e,e^2,e^3,\ldots] ={}&  1e^0 + 1e^1 + 2e^2 + 3e^3 + 5e^4 + 7e^5 + 11e^6  
\nonumber \\ & + 15e^7 + 22e^8 + 30e^9 + 42e^{10} + ... 
=
\sum_{n \ge 0} p(n)  \, e^n
\end{align}

With all coefficients unity, the edge layering of $\Mb$ does not appear to have the generality or usefulness of Fine's lemma.

\section{Conclusion}

The combinatorics of types of tubdigons tells us their generating series is the zero of a polynomial.
Solving that polynomial with Wildberger's polynomial formula and applying the multinomial theorem allows us to equate terms, recovering Fine's Identity.

Ira Gessel (personal communication) tells me that this sort of Fine Identity proof is well known.

\bibliographystyle{vancouver}
\bibliography{HyperCatBib.bib}

\begin{bibdiv}
\begin{biblist}

\bib{Erdelyi1940}{article}{
      author={Erd{\'e}lyi, Arthur},
      author={Etherington, I. M.~H.},
       title={Some problems of non-associative combinations (2)},
        date={1940},
     journal={Edinburgh math notes},
      volume={32},
       pages={vii\ndash xii},
}

\bib{Fine1959}{article}{
      author={Fine, Nathan~J.},
       title={Sums over partitions},
        date={1959},
        ISSN={issn},
     journal={Report of the Institute in the Theory of Numbers},
      volume={June 21–July 17},
       pages={86–94},
        note={Reproduced as chapter 2 in \cite{Fine1988}},
}

\bib{Fine1988}{book}{
      author={Fine, Nathan~J.},
       title={Basic hypergeometric series and applications},
      series={Mathematical Surveys and Monographs},
   publisher={American Mathematical Society},
     address={Providence, RI},
        date={1988},
      volume={27},
        ISBN={0-8218-1524-5},
      review={\MR{956465 (91j:33011)}},
}

\bib{Kreweras1972}{article}{
      author={Kreweras, G.},
       title={Sur les partitions non croisees d'un cycle [on the non-crossing partitions of a cycle]},
        date={1972},
        ISSN={0012-365X},
     journal={Discrete math},
      volume={1},
      number={4},
       pages={333\ndash 350},
}

\bib{KrewerasEarnshaw2005}{misc}{
      author={Kreweras, Germain},
      author={Earnshaw, Berton~A.},
       title={On the non-crossing partitions of a cycle},
        date={2005},
         url={https://users.math.msu.edu/users/earnshaw/research/kreweras.pdf},
        note={English translation of Kreweras~\cite{Kreweras1972}},
}

\bib{Raney1960}{article}{
      author={Raney, George~N.},
       title={Functional composition patterns and power series reversion},
        date={1960},
     journal={Trans am math soc},
      volume={94},
       pages={441\ndash 451},
}

\bib{Rubine2025Finite}{misc}{
      author={Rubine, Dean},
      author={Mukewar, Pratham},
       title={Finite interpretation of series zeros},
        date={2025},
        note={forthcoming},
}

\bib{Tutte1964}{article}{
      author={Tutte, W.~T.},
       title={The number of planted plane trees with a given partition},
        date={1964},
        ISSN={00029890, 19300972},
     journal={Am math mon},
      volume={71},
      number={3},
       pages={272\ndash 277},
}

\bib{Wildberger2025}{article}{
      author={Wildberger, N.~J.},
      author={Rubine, Dean},
       title={A hyper-catalan series solution to polynomial equations, and the geode},
        date={2025},
     journal={The American Mathematical Monthly},
      volume={132},
      number={5},
       pages={383\ndash 402},
      eprint={https://doi.org/10.1080/00029890.2025.2460966},
         url={https://doi.org/10.1080/00029890.2025.2460966},
}

\end{biblist}
\end{bibdiv}
\end{document}